\documentclass[10pt]{amsart}

\usepackage{tikz}
\usetikzlibrary{positioning}
\usetikzlibrary{arrows}
\usepackage{amssymb,amsmath,latexsym,graphicx}
\usepackage{charter,eucal}
\usepackage{amssymb}
\usepackage{amscd}
\usepackage{tcolorbox}
\usepackage[all,cmtip]{xy}
\usepackage{rotating,multicol,soul}
\usetikzlibrary{decorations.markings}

\setul{}{1.1pt}

\newtheorem{theorem}{Theorem }[section]

\newtheorem{lemma}[theorem]{Lemma}
\newtheorem{observation}[theorem]{Observation}

\newtheorem{remark}[theorem]{Remark}
\newtheorem{corollary}[theorem]{Corollary}
\newtheorem{proposition}[theorem]{Proposition}
\newtheorem{principle}[theorem]{\textsc{Principle}}

\newcommand{\bt}{\begin{theorem}}
\newcommand{\et}{\end{theorem}}
\newcommand{\bmt}{\begin{maintheorem}}
\newcommand{\emt}{\end{maintheorem}}
\newcommand{\bc}{\begin{corollary}}
\newcommand{\bl}{\begin{lemma}}
\newcommand{\ec}{\end{corollary}}
\newcommand{\el}{\end{lemma}}
\newcommand{\bo}{\begin{observation}}
\newcommand{\eo}{\end{observation}}
\newcommand{\bp}{\begin{proposition}}
\newcommand{\ep}{\end{proposition}}
\newcommand{\br}{\begin{remark}}
\newcommand{\er}{\end{remark}}
\newcommand{\bpr}{\begin{principle}}
\newcommand{\epr}{\end{principle}}

\def\I{\mathop{\mathrm{I}}}

\def\I{\mathbf{I}}

\def\eop{\hspace*{\fill}$\blacksquare$}

\newcommand{\F}{\mathbb{F}}
\newcommand{\bP}{\mathbb{P}}

\newcommand{\mG}{\mathcal{G}}

\newcommand{\mQ}{\mathcal{Q}}

\newcommand{\mP}{\mathcal{P}}
\newcommand{\mL}{\mathcal{L}}
\newcommand{\mE}{\mathcal{E}}

\newcommand{\mW}{\mathcal{W}}

\newcommand{\mO}{\mathcal{O}}

\newcommand{\mS}{\mathcal{S}}
\newcommand{\mD}{\mathcal{D}}

\newcommand{\cV}{\mathcal{V}}

\newcommand{\mB}{\mathcal{B}}
\newcommand{\mH}{\mathcal{H}}

\newcommand{\hT}{\mathbf{T}}

\newcommand{\wt}{\widetilde}
\newcommand{\ol}{\overline}

\usetikzlibrary{matrix}
\usetikzlibrary{arrows}
\usepackage{pgf}
\usepackage{lscape}
\usepackage{listings}

\usepackage{enumitem}

\title{Epimorphisms of generalized polygons A: the planes, quadrangles and hexagons}
\keywords{generalized polygon, epimorphism, cover}
\author{Joseph A. Thas and Koen Thas}
\thanks{}
\address{Ghent University, Department of Mathematics, Krijgslaan 281, S25, B-9000, Ghent, Belgium}
\email{joseph.thas@ugent.be; koen.thas@gmail.com}\date{}

\begin{document}

\maketitle
\begin{abstract}
Inspired by a theorem by Skornjakov--Hughes--Pasini \cite{Skorn,H,P} and a problem which turned up in our recent paper \cite{part2}, we start a study of epimorphisms with source a thick generalized $m$-gon and target a thin generalized $m$-gon. In this first part of the series, we classify the cases $m = 3, 4$ and $6$ when the polygons are finite. Then we show that the infinite case is very different, and construct examples which strongly deviate from the finite case. A number of general structure theorems are also obtained. We introduce the theory of locally finitely generated generalized polygons and locally finitely chained generalized polygons along the way.\\
\end{abstract}

\begin{tcolorbox}
\tableofcontents
\end{tcolorbox}

\section{Introduction}

The notion of epimorphism has played a fundamental role in many areas of Mathematics for a long time now. In this paper, we want to study epimorphisms in the particular context of generalized polygons. Our main motivation originates from an interesting result of Pasini \cite{P}, which generalizes an older result of Skornjakov \cite{Skorn} and Hughes \cite{H} (the latter two handle the case of projective planes). Here is the statement. 

\begin{theorem}[Skornjakov--Hughes--Pasini \cite{Skorn,H,P}]
\label{HP}
Let $\alpha$ be a morphism from a thick (possibly infinite) generalized $m$-gon $\mE$ to a thick (possibly infinite) generalized $m$-gon $\mE'$, with $m \geq 3$. If $\alpha$ is surjective, then 
either $\alpha$ is an isomorphism, or each element in $\mE'$ has an infinite fiber in $\mE$.
\end{theorem}

(We will formally define a number of notions in the next section.)
As a direct corollary, any epimorphism between finite thick generalized polygons must be an isomorphism. Of course, from the viewpoint of the category of finite sets this result might come with little surprise, but the extra geometrical structure makes such results less trivial. As Van Maldeghem remarks in \cite[section 4.2.4]{POL}, the thickness assumption is crucial here, and in that same remark he briefly mentions a counter example in the thin case. Wondering what the ``thin version'' of Theorem \ref{HP} is, is one of our main starting points.  

A local variation on Theorem \ref{HP} by B\"{o}di and Kramer \cite{BK} states that an epimorphism between thick generalized 
$m$-gons ($m \geq 3$) is an isomorphism if and only if its restriction to at least one point row or line pencil is bijective. Later on, Gramlich and Van Maldeghem thoroughly studied epimorphisms from thick generalized $m$-gons to thick generalized $n$-gons with $n < m$ in their works \cite{GramHVM0,GramHVM}, and again, classification results are obtained based on the local nature of the epimorphisms.  
Finally, we mention recent work of Koen Struyve \cite{KS} on epimorphisms of spherical Moufang buildings. One of the main results of that paper is that the epimorphic image of a Moufang generalized polygon necessarily again is a Moufang polygon. \\

A second motivation came from our recent paper \cite{part2}, which is a second instalment in a series of papers \cite{part1,part2} about covers of generalized quadrangles. In that study, a question naturally popped up involving morphisms between thick finite generalized quadrangles which surjectively map some geometrical hyperplane of the source quadrangle to a geometrical hyperplane of the target quadrangle. The precise statement is described in the second last section of this paper.  

\subsection*{Synopsis of the present paper}

In section \ref{basic}, we introduce some notions which will be frequently used throughout the paper. In sections \ref{plan}|\ref{quad}|\ref{hexa}, we classify epimorphisms from thick finite generalized $m$-gons with $m \in \{ 3, 4, 6 \}$ to thin generalized $m$-gons with an order. 

In section \ref{chain} we introduce locally finitely generated generalized polygons and locally finitely chained generalized polygons and provide a number of examples which in their turn provide classes of epimorphisms which agree with  our results in the finite case. On the other hand,  in section \ref{counter} we construct classes of examples of epimorphisms which deviate from our results in the finite case.   
In section \ref{lf} we explain corollaries of the proofs of the results in section \ref{quad} and section \ref{hexa} for locally finite quadrangles and hexagons. Then, in section \ref{thin}, we study the particular case of epimorphisms between thin generalized polygons. Finally, in section \ref{app}, we solve the question mentioned in our paper \cite{part2}.

\section{Some basic definitions}
\label{basic}

We summarize a number of definitions which we will need in due course. 

\subsection{Generalized polygons}

Let $\Gamma = (\mP,\mL,\I)$ be a point-line geometry, and let $m$ be a positive integer at least $2$. We say that $\Gamma$ is a {\em weak generalized $m$-gon} if any two elements in $\mP \cup \mL$ are contained in at least one ordinary sub $m$-gon (as a subgeometry of $\Gamma$), and if $\Gamma$ does not contain ordinary sub $k$-gons with $2 \leq k < m$. For $m = 2$ every point is incident with every line. 

If $m \geq 3$, we say $\Gamma$ is a {\em generalized $m$-gon} if furthermore $\Gamma$ contains an ordinary sub $(m + 1)$-gon as a subgeometry. Equivalently, a weak generalized $m$-gon with $m \geq 3$ is a generalized $m$-gon if it is {\em thick}, meaning that every point is incident with at least three distinct lines and every line is incident with at least three distinct points. A weak generalized $m$-gon is {\em thin} if it is not thick; in that case, we also speak of {\em thin generalized $m$-gons}. If we do not specify $m$ (the ``gonality''), we speak of {\em (weak) generalized polygons}. Note that thick generalized $2$-gons (or {\em generalized digons}) do not contain ordinary $3$-gons as a subgeometry.  
 
It can be shown that generalized polygons have an {\em order} $(u,v)$: there exists positive integers $u \geq 2$ and $v \geq 2$ such that each 
point is incident with $v + 1$ lines and each line is incident with $u + 1$ points. We say that a weak generalized polygon is {\em finite} if its number of points and lines is finite | otherwise it is {\em infinite}.  If a thin weak generalized polygon has an order $(1,u)$ or $(u,1)$ it is called a {\em thin 
generalized polygon} of order $(1,u)$ or $(u,1)$. 

Note that the generalized $3$-gons are precisely the (axiomatic) projective planes. Generalized $4$-gons, resp. $6$-gons, resp. $8$-gons are also called {\em generalized quadrangles}, resp. {\em hexagons}, resp. {\em octagons}.

\subsection{Sup polygons}

We say that $\Gamma' = (\mP',\mL',\I')$ is a {\em sub generalized $m$-gon} of the generalized $m$-gon $\Gamma$, $m \geq 3$, if $\mP' \subseteq \mP$, $\mL' \subseteq \mL$, and if $\I'$ is the induced incidence coming from $\I$.

\subsection{Morphisms and epimorphisms}

A {\em morphism} from a weak generalized polygon $\Gamma = (\mP,\mL,\I)$ to a weak generalized polygon $\Gamma = (\mP',\mL',\I')$ is a map $\alpha: \mP \cup \mL \mapsto \mP' \cup \mL'$ which maps points to points, lines to lines and which preserves the incidence relation (note that we do not ask the gonalities to be the same). We say that a morphism $\alpha$ is an {\em epimorphism} if $\alpha(\mP) = \mP'$ and $\alpha(\mL) = \mL'$. Contrary to Gramlich and Van Maldeghem \cite{GramHVM0,GramHVM}, we do not ask surjectivity onto the set of flags of $\Gamma'$ (the incident point-line pairs of $\Gamma'$). 

If an epimorphism is injective, and if the inverse map is also a morphism, then we call it an {\em isomorphism}.

\subsection{Doubling}

Let $\Gamma = (\mP, \mL, \I)$ be a (not necessarily finite) generalized $n$-gon of order $(s,s)$ (for $n = 3$, projective planes of order $(1,1)$ are allowed). Define the {\em double of $\Gamma$} as the generalized $2n$-gon $\Gamma^\Delta$ which arises by letting its point set be $\mP \cup \mL$, and letting its line set be the flag set of $\Gamma$ (the set of incident point-line pairs).  Its parameters are $(1, s)$. The full automorphism group of $\Gamma^\Delta$ is isomorphic to the  group consisting of all automorphisms and dualities (anti-automorphisms) of $\Gamma$. Sometimes we prefer to work in the point-line dual of $\Gamma^\Delta$, but we use the same notation (while making it clear in what setting we work). This is what we will do in this section.
Vice versa, if $\Gamma'$ is a thin generalized $2n$-gon of order $(1,s)$, then it is isomorphic to the double $\Gamma^\Delta$ of a generalized $n$-gon $\Gamma$ of order $(s,s)$. \\

\medskip
\section{Epimorphisms to thin planes}
\label{plan}

In this first research section, we classify epimorphisms with source a thick projective plane and target a thin projective plane of order $(1,1)$.

\begin{theorem}
\label{GT}
Let $\Phi$ be an epimorphism of a thick projective plane $\mP$ onto a thin projective plane $\Delta$ of order $(1,1)$.  Then exactly two classes of epimorphisms $\Phi$ occur (up to a suitable permutation of the points of $\Delta$), and they are described as follows. 

\begin{itemize}
\item[{\rm(a)}] 
The points of $\Delta$ are $\overline{a}, \overline{b}, \overline{c}$, with $\overline{a} \sim \overline{b} \sim \overline{c} \sim \overline{a}$, and put $\Phi^{-1} (\overline{x}) = \wt X$, with $\overline{x} \in \lbrace \overline{a}, \overline{b}, \overline{c} \rbrace$.

   Let $(\wt{A}, \wt{B})$, with $\wt{A} \ne \emptyset \ne \wt{B}$, be a partition of the set of all points incident with a line $L$ of $\mP$. Let  $\wt{C}$ consist of the points not incident with $L$. Furthermore, $\Phi^{-1}(\overline{a}\overline{b}) = L$, $\Phi^{-1}(\overline{b}\overline{c})$ is the set of all lines distinct from $L$ but incident with a point of $\wt{B}$ and $\Phi^{-1}(\overline{a}\overline{c})$ is the set of all lines distinct from $L$ but incident with a point of $\wt{A}$.
   \item[{\rm(b)}] 
   The dual of (a).
\end{itemize}
\end{theorem}

{\em Proof}.\quad
We proceed in a number of steps. Below, if $U$ is any line of $\mP$,  $\widehat{U}$ is the set of all points incident with $L$. \\

\quad (1) {\bf Let $L$ be any line of $\mP$. Then $\Phi$ maps the set of all points incident with $L$ onto the set of all points incident with $\Phi(L)$.}

{\em Proof of {\rm (1)}}.\quad
Suppose w. l. o. g. that $\Phi(L) = \overline{a}\overline{b}$. 
Assume that $\Phi(\widehat{L}) = \{\overline{a}\}$. Let $\Phi(x) = \overline{b}$. Then if $y \I L$, it follows that $\Phi(\widehat{xy}) = \{ \overline{a}, \overline{b}\}$, so that 
all points of $\mP$ are mapped to a point in $\{ \overline{a}, \overline{b} \}$, contradiction.    \eop \\

Of course, the dual property also holds.\\

\quad (2) {\bf For at least one point $u$ of $\mP$ it is not possible that  at least two lines 
of $\mP$ incident with $u$ are mapped onto $\overline{u}\overline{v}$ and at least two lines of $\mP$ incident with $u$ are mapped onto $\overline{u}\overline{w}$, with $\overline{u} = \Phi(u)$ and $\{ \overline{u}, \overline{v}, \overline{w} \} = \{ \overline{a}, \overline{b}, \overline{c} \}$. } \\

{\em Proof of {\rm (2)}}.\quad 
Assume the contrary. 

Take $L, M \in \Phi^{-1}(\overline{a}\overline{b})$ with $L \ne M$. 
\begin{itemize}
\item
Suppose that $L \cap M = \{ w\}$ with $w \in \widetilde{A}$. Let $w' \in \widetilde{C}$. Then if $r \I L$ with $r \in \widetilde{A}$, we have that $r' \in \widetilde{A}$ where $\{ r' \} = M \cap w'r$. 
It follows that $\vert \widehat{L} \cap \widetilde{A} \vert = \vert \widehat{M} \cap \widetilde{A} \vert$. 
\item
Now suppose $w \in \widetilde{B}$. Let $u' \in \widetilde{C}$. Then if $r \I L$ with $r \in \widetilde{A}$, it follows that $r' \in \widetilde{A}$ where $\{r' \} = M \cap w'r$. So $\vert \widehat{L} \cap \widetilde{A} \vert = \vert \widehat{M} \cap \widetilde{A} \vert$. 
\end{itemize}

Now take $L' \in \Phi^{-1}(\overline{a}\overline{b})$ and $M' \in \Phi^{-1}(\overline{a}\overline{c})$. Then $L ' \cap M' = \{w\} \subset \widetilde{A}$. Now let $w' \in 
\widetilde{C}$ with $w'$ not incident with $M'$ (by assumption there is a line $M'' \ne M'$ in $\mP$ incident with $w$ and $\Phi(M'') = \overline{a}\overline{c}$, and 
on $M''$ there is a point $w'$ mapped onto $\overline{c}$).  
Then if $r' \in \widehat{M'} \cap \widetilde{A}$, it follows that $r \in \widetilde{A}$ where $\{ r \} = L' \cap w'r'$. So $\vert \widehat{L'} \cap \widetilde{A} \vert \geq \vert \widehat{M'} \cap \widetilde{A} \vert$. Doing the same reasoning with $w'' \in \widetilde{B}$ instead of $w' \in \widetilde{C}$, we obtain that $\vert \widehat{L'} \cap \widetilde{A} \vert \leq \vert \widehat{M'} \cap \widetilde{A} \vert$ and hence equality holds.

We have shown that every line in $\mP$ meets $\widetilde{A}$ in $0$ or $\alpha$ (a constant) points. Similarly we have that every line in $\mP$ meets $\widetilde{B}$ in $0$ or $\beta$ (a constant) points and that every line in $\mP$ meets $\widetilde{C}$ in $0$ or $\gamma$ (a constant) points.

If $\mP$ has order $n$, we now have the following equalities:

\begin{equation*}
\left\{
\begin{aligned}
  \alpha + \beta&= n + 1\\
  \alpha + \gamma &= n + 1\\
  \beta + \gamma &= n + 1
\end{aligned}
\right.
\end{equation*}

It follows that $\alpha = \beta = \gamma = \frac{n + 1}{2}$. 

Now take any point $a \in \widetilde{A}$ and count the points on $\widetilde{A}$ by considering the points of $\widetilde{A}$ on the lines incident with $a$ to obtain that $\vert \widetilde{A} \vert 
= (n + 1)(\alpha - 1) + 1 = \frac{n^2 + 1}{2}$. Next, let $a'$ be a point of $\mP$ which is not contained in $\widetilde{A}$; then considering the intersections of 
the lines incident with $a'$ and $\widetilde{A}$, we also see that $\frac{n + 1}{2}$ divides $\vert \widetilde{A} \vert$, contradiction. \\

\quad (3) {\bf Ending the proof.}\\
By (2) we may suppose that $x \in \widetilde{A}$ is incident with just one line $L$ of $\Phi^{-1}(\overline{a}\overline{b})$. It follows that $\widetilde{B} \subset \widehat{L}$. 
We distinguish two cases.
\begin{itemize}
\item[(a)]
Suppose $\widetilde{A} \subset \widehat{L}$. Then it is easy to see that we are in Case (a) of the theorem.
\item[(b)]
Suppose $\widetilde{A} \not\subset \widehat{L}$.
Let $y \in \widetilde{A}$ be not incident with $L$. First suppose $\vert \widetilde{B} \vert \geq 2$, and let $b, b'$ be 
different points of $\widetilde{B}$. Then on $by$ we have $n$ points of $\widetilde{A}$. By joining $b'$ with all these points, we see that all points of $\mP$ which 
are not incident with $L$ are contained in $\widetilde{A}$. But that is a contradiction since we now would not have points in $\widetilde{C}$. It follows 
that $\vert \widetilde{B} \vert = 1$, which implies that we are in Case (b) of the theorem.
\end{itemize}

\eop \\

\medskip
\section{Epimorphisms to thin quadrangles}
\label{quad}

In this section, we classify epimorphisms with source a thick generalized quadrangle and target a thin generalized quadrangle of order $(s',1)$.

\begin{theorem}
\label{JATGQ}
Let $\Phi$ be an epimorphism of a thick generalized quadrangle $\mathcal{S}$ of order $(s, t)$ onto a grid $\mathcal{G}$. Let $\mathcal{G}$ have order $(s^\prime, 1)$. Then $s^\prime = 1$ and exactly two classes of epimorphisms $\Phi$ occur (up to a suitable permutation of the points of $\mathcal{G})$.

\begin{itemize}
\item[{\rm(a)}] The points of $\mathcal{G}$ are $\overline{a}, \overline{b}, \overline{c}, \overline{d}$, with $\overline{a} \sim \overline{b} \sim \overline{c} \sim \overline{d} \sim \overline{a} $, and put $\Phi^{-1} (\overline{x}) = \wt X$, with $\overline{x} \in \lbrace \overline{a}, \overline{b}, \overline{c}, \overline{d} \rbrace$.

   Let $(\wt{A}, \wt{B})$, with $1 \le \vert \wt A \vert \le s, 1 \le \vert \wt B \vert \le s$, be a partition of the set of all points incident with a line $L$ of $\mathcal{S}$. Let  $\wt{C}$ consist of the points not incident with $L$ but collinear with a point of $\wt{B}$, and let $\wt{D}$ consist of the points not incident with $L$ but collinear with a point of $\wt{A}$. Further, $\Phi^{-1}(\overline{a}\overline{b}) = \{ L \}$, $\Phi^{-1}(\overline{b}\overline{c})$ is the set of all lines distinct from $L$ but incident with a point of $\wt{B}$, $\Phi^{-1}(\overline{a}\overline{d})$ is the set of all lines distinct from $L$ but incident with a point of $\wt{A}$ and $\Phi^{-1}(\overline{c}\overline{d})$ consists of all lines incident with at least one point of $\wt{C}$ and at least one point of $\wt{D}$.
\item[{\rm(b)}] The dual of (a).
\end{itemize}
\end{theorem}

{\it Proof}.\quad We proceed in a number of steps.

\quad (\rm 1)  {\it {\bf Let $L$ be any line of $\mathcal{S}$. Then $\Phi$ maps the set of all points incident with $L$ onto the set of all points incident with $\Phi(L)$}.}

{\em Proof of {\rm (1)}}.\quad 
Let $L$ be any line of $\mathcal{S}$, let $\Phi(L) = \overline{L}$, let $\widehat{L}$ be the set of all points incident with $L$, and let $ \widehat {\overline{L}}$ be the set of all points incident with $\overline{L}$. Assume, by way of contradiction, that $\Phi(\widehat{L}) \not = \widehat{\overline{L}}$. Let $\overline{x}  \:  \mathbf{I} \:  \overline{L}, \overline{x} \notin \Phi(\widehat{L})$. Let $\overline{x} \sim \overline{y}, \overline{y} \not\mathbf{I} \: \overline{L}$ and $y \in \Phi^{-1}(\overline{y})$. If $z \: \mathbf{I} \:  L, z \sim y$, then $\Phi(z) = \overline{x}$. So $\overline{x} \in \Phi(\widehat{L})$, a contradiction.\\

\quad (\rm 2) {\it {\bf  Let $x$ be any point of $\mathcal{S}$, and let $\Phi(x) = \overline{x}$. Then the image of the set of lines incident with $x$, is the set of lines incident with $\overline{x}$}.}

{\em Proof of {\em (2)}}\quad Dualize the proof of (1).\\

\quad (\rm 3) {\bf We have that $s^\prime = 1$.}

{\em Proof of {\rm (3)}}.\quad Let $A, B$ be the two systems of lines of $\mathcal{G}$, let $\Phi^{-1}(A)$ be the union of all inverse images of the elements of $A$, and let $\Phi^{-1}(B)$ be the union of all inverse images of the elements of $B$. If $x$ is a point of $\mathcal{S}$, then $A_x$ is the set of all lines of $\Phi^{-1}(A)$ incident with $x$ and $B_x$ is the set of all lines of $\Phi^{-1}(B)$ incident with $x$.
   
   Let $c$ be a point of $\mathcal{S}$ and let $M$ be a line incident with $c$. Let $\overline{c} = \Phi(c)$, let $\overline{d}$  be a point of $\mathcal{G}$ not collinear with $\overline{c}$, and let $\Phi(d) = \overline{d}$. Then $c \not \sim d$. The line $\overline{M} = \Phi(M)$ is incident with $\overline{c}$; say $\overline{M} \in B$. Let $\overline{d} \: \mathbf{I} \: \overline{U} \sim \overline{M}$ and $d \: \mathbf{I}\: U \sim M$. Then $U \in \Phi^{-1}(\overline{U})$ and $\overline{U} \in A$. It easily follows that $\vert B_c \vert = \vert A_d \vert$.
   
   Assume, by way of contradiction, that $s^\prime > 1$. Hence there is a point $\overline{e}$ of $\mathcal{G}$ with $\overline{e} \not \sim \overline{c}$ and $\overline{e} \not \sim \overline{d}$. So $e \not \sim c$ and $e \not \sim d$, with $e \in \Phi^{-1}(\overline{e})$. Consequently $\vert B_c \vert = \vert A_d \vert = \vert B_e \vert = \vert A_c \vert$.
   
   Hence $\vert A_c \vert = \vert B_c \vert$ for each point $c$ of $\mathcal{S}$. It follows that $t$ is odd, so $t \ge 3$.
   
   Let $r$ be a point of $\mathcal{S}$. Then $\vert A_r \vert = \vert B_r \vert = \frac{t+1}{2} \ge 2$. Let $L_1, L_2 \in \Phi^{-1}(A), r \: \mathbf{I} \: L_1, r \: \mathbf{I} \: L_2, L_1 \not = L_2$. Further, let $u_1 \: \mathbf{I} \: L_1, u_2 \: \mathbf{I} \: L_2, u_1 \not = r \not = u_2, \Phi(L_1) = \Phi(L_2) = \overline{L}, \Phi(r) = \overline{r} \: \mathbf{I} \: \overline{L}, \Phi(u_1) = \overline{u}_1 \: \mathbf{I} \: \overline{L}, \overline{u}_1 \not = \overline{r}, \Phi(u_2) = \overline{u}_2 \: \mathbf{I} \overline{L} , \overline{u}_1 \not = \overline{u}_2$. Assume that $u_1 \: \mathbf{I} \: U_1, U_1 \in B_{u_1}, u_2 \: \mathbf{I} \: U_2, U_1 \sim U_2$. Then $\Phi(U_2) = \overline{L}  \in A$. Hence $\vert B_{u_1} \vert \leq \vert A_{u_2} \setminus \lbrace L_2 \rbrace\vert$. So $\vert B_{u_1} \vert \ < \vert A_{u_2} \vert$, a contradiction as $\vert B_{u_1} \vert = \vert A_{u_2} \vert = \frac{t+1}{2}$.
   
   We conclude that $s^\prime = 1$.\\
   
\quad{\rm (4)} {\it {\bf Let $s^\prime = 1$. Then not for any given points $\overline{x}$ and $\overline{y}$ in $\mathcal{G}$, with $\overline{x} \sim \overline{y}, \overline{x} \not = \overline{y}$, there exist points $x \in \Phi^{-1}(\overline{x})$ and $y \in \Phi^{-1}(\overline{y})$ with $x \not \sim y$}.}
   
Let $\overline{u} \not \sim \overline{v}$, with $\overline{u} = \Phi(u)$ and $\overline{v} = \Phi(v)$. Then by (3) we have $\vert B_u \vert = \vert A_v \vert$ and $\vert A_u \vert = \vert B_v \vert$.

{\em Proof of {\rm (4)}}.\quad Assume, by way of contradiction, that for any given points $\overline{x}$ and $\overline{y}$ in $\mathcal{G}$, with $\overline{x} \sim \overline{y}, \overline{x} \not = \overline{y}$, there exist points $x \in \Phi^{-1}(\overline{x})$ and $y \in \Phi^{-1}(\overline{y})$ with $x \not \sim y$.
   
   Let $\overline{a}, \overline{b}, \overline{c}, \overline{d}$ be the points of $\mathcal{G}$, with $\overline{a} \sim \overline{b} \sim \overline{c} \sim \overline{d} \sim \overline{a}$, and let $\overline{a} = \Phi(a)$. By (1) and (2) we can choose points $d, c$ such that $a \sim d \sim c$, with $\overline{d} = \Phi(d)$ and $\overline{c} = \Phi(c)$. Similarly we can choose a line $M$ incident with $c$ such that $\Phi(M) = \ol{b}\ol{c}$. Let $b \: \mathbf{I} \: M, a \sim b$. Then clearly $\Phi(b) = \ol{b}$. So $a \sim b \sim c \sim d \sim a$. Let $A = \lbrace \ol{b}\ol{c}, \ol{a}\ol{d} \rbrace$ and $B = \lbrace \ol{a}\ol{b}, \ol{c}\ol{d} \rbrace$.
   
   It follows that $\vert A_c \vert = \vert B_a \vert, \vert B_c \vert = \vert A_a \vert, \vert A_b \vert = \vert B_d \vert, \vert B_b  \vert = \vert A_d \vert$.
   
   Let $\Phi(a^\prime) = \ol{a}$. As $\ol{a} \not \sim \ol{c}$, we have $a^\prime \not \sim c$. Hence $\vert A_a \vert = \vert B_c \vert = \vert A_{a^\prime} \vert$ and $\vert B_a \vert = \vert A_c \vert = \vert B_{a^\prime} \vert $.
   
   By assumption there is a point $b^\prime \not \sim a$ with $\Phi(b^\prime) = \ol{b}$. Let $b^\prime \sim r \sim a$, with $\Phi(ra) \in A$. Then $\Phi(r) = \ol{a}$, so $\Phi(rb^\prime) = \ol{a}\ol{b} \in B$. Hence $\vert B_{b^\prime} \vert \ge \vert A_a \vert$, so $\vert B_b \vert \ge \vert A_a \vert$ and thus $\vert A_d \vert \ge \vert A_a \vert$. Similarly we obtain $\vert A_a \vert \ge \vert A_d \vert$. It follows that $\vert A_a \vert = \vert A_d \vert$.
   
   Now clearly $\vert A_a \vert = \vert A_b \vert =\vert A_c \vert =\vert A_d \vert$. Similarly $\vert B_a \vert = \vert B_b \vert = \vert B_c \vert = \vert B_d \vert$. Consequently $\vert A_u \vert = \vert B_v \vert$ for all $u, v \in \lbrace a, b, c, d \rbrace$ and $\vert A_w \vert = \vert B_w \vert = \frac{t+1}{2}$ for all $w \in \lbrace a, b, c, d \rbrace$. Hence $\vert A_w \vert = \vert B_w \vert = \frac{t+1}{2}$ for each point $w$ of $\mS$.
   
   Such as in the last part of (3) we now find a contradiction.\\
   
\quad{\rm (5}) {\it {\bf For $s^\prime = 1$ two (mutually dual) cases occur}.}
   
{\em Proof of {\rm (5)}}.\quad   By (4) we know that not for any given points $\ol{x}$ and $\ol{y}$ in $\mG$, with $\ol{x} \sim \ol{y}, \ol{x} \not = \ol{y},$ there exist points $x \in \Phi^{-1}(\ol{x})$ and $y \in \Phi^{-1}(\ol{y})$ with $x \not \sim y$.
   
   Let $\ol{a}, \ol{b}, \ol{c}, \ol{d}$ be the points of $\mG$, with $\ol{a} \sim \ol{b} \sim \ol{c} \sim \ol{d} \sim \ol{a}$, let $A = \lbrace \ol{b}\ol{c}, \ol{a}\ol{d} \rbrace$ and $B = \lbrace \ol{a}\ol{b}, \ol{c}\ol{d} \rbrace$, and let $\Phi^{-1}(\ol{x}) = \wt{X}$ for $\ol{x} \in \lbrace \ol{a}, \ol{b}, \ol{c}, \ol{d} \rbrace$.
   
   Without loss of generality, we assume that there do not exist points $a^\prime \in \wt{A}$ and $b^\prime \in \wt{B}$ for which $a^\prime \not \sim b^\prime$. Let $a^\ast \in \wt{A}, b^\ast \in \wt{B}$. If $e \I a^\ast b^\ast$, then $\Phi(e) \in \lbrace \ol{a}, \ol{b} \rbrace$. Hence each point of $a^\ast b^\ast$ belongs to $\wt{A} \cup \wt{B}$.
   
   \begin{itemize}
   \item[(5.a)]{\it Assume that there is a line $a^\ast b^\ast$ which has at least two points in common with $\wt{A}$ and at least two points in common with $\wt{B}$}.
   
   Assume, by way of contradiction, that $\wt{A} \cup \wt{B}$ contains a point $g$ not incident with the line  $a^\ast b^\ast$. Let $g \in \wt{A}$. Then $g$ is collinear with at least two points of $a^\ast b^\ast$ in $\wt{B}$, yielding a triangle, clearly a contradiction. Consequently, $\wt{A} \cup \wt{B}$ is the set of all points of $\mS$ incident with the line $a^\ast b^\ast$. Then $\wt{C}$ is the set of all points of $\mS$ not incident with $a^\ast b^\ast$ but collinear with a point of $\wt{B}$, and $\wt{D}$ is the set of all points of $\mS$ not incident with $a^\ast b^\ast$ but collinear with a point of $\wt{A}$. Hence $\Phi^{-1}(\ol{a}\ol{b}) = \lbrace L \rbrace, \Phi^{-1}(\ol{b} \ol{c})$ is the set of all lines distinct from $L$ but incident with a point of $\wt{B}$, $\Phi^{-1}(\ol{a}\ol{d})$ is the set of all lines distinct from $L$ but incident with a point of $\wt{A}$ and $\Phi^{-1}(\ol{c}\ol{d})$ consists of all lines incident with at least one point of $\wt{C}$ and at least one point of $\wt{D}$.
   
   Hence we have Case (a) in the statement of the theorem.
   \item[(5.b)]{\it Assume that each line $a^\ast b^\ast$ is incident with exactly one point of one of the sets $\wt{A}, \wt{B}$, and with $s$ points of the other set}.
   
   Let $a^\ast \in \wt{A}, b^\ast \in \wt{B}$ and assume that $a^\ast$ is the only point of $\wt{A}$ which is incident with the line $a^\ast b^\ast$. If $e \in \wt{A} \setminus \lbrace a^\ast \rbrace$, then $e$ is collinear with at least $s$ points incident with the line $a^\ast b^\ast$, clearly a contradiction. Hence $\wt{A} = \lbrace a^\ast \rbrace$. Then $\wt{B}$ consists of all points, distinct from $a^\ast$, incident with the lines of a set $\cV$ consisting of $u \: (\ge 1)$ lines incident with $a^\ast$. 
   
   Let $\mW$ be the set of the $t + 1 - u$ lines which do not belong to $\mathcal{V}$ and which are incident with $a^\ast$. It is clear that the points, distinct from $a^\ast$, which are incident with these $t + 1 - u$ lines are mapped by $\Phi$ onto $\ol{d}$. Assume, by way of contradiction, that $h \not \sim a^\ast$ is mapped onto $\ol{d}$. Let $h^\prime$ be collinear with $h$ and incident with some line of $\cV$. Then $\ol{d} = \Phi(h) \sim \Phi(h^\prime) = \ol{b}$, a contradiction. Hence $\wt{D}$ consists of all points, distinct from $a^\ast$, incident with the $t + 1 - u$ lines of $\mW$. It follows that $u \le t$. Finally, $\wt{C} = \lbrace c^\ast \Vert  \: c^\ast \not \sim a^\ast \rbrace$.
  
   Consequently we have the dual of Case (a), that is, Case (b) in the statement of the theorem. \eop \\
   \end{itemize}

 \section{Epimorphisms to thin hexagons}  
 \label{hexa}

  In the next theorem, we classify epimorphisms from finite thick generalized hexagons to thin generalized hexagons of order $(s',1)$. As in the case of Theorem \ref{JATGQ}, a part of this theorem can be taken over in the infinite case.  
   
   \begin{theorem}
   \label{JATGH}
   Let $\Phi$ be an epimorphism of a thick generalized hexagon $\mS$ of order $(s, t)$ onto a thin generalized hexagon $\mG$ of order $(s^\prime, 1)$. Then $s^\prime = 1$ and exactly two classes of epimorphisms $\Phi$ occur (up to a suitable permutation of the points of $\mG$).
   
   \begin{itemize}
   \item[{\rm (a)}] The points of $\mG$ are $\ol{a}, \ol{b}, \ol{c}, \ol{d}, \ol{e}, \ol{f}$, with $\ol{a} \sim \ol{b} \sim \ol{c} \sim \ol{d} \sim \ol{e} \sim \ol{f} \sim \ol{a}$, and put $\Phi^{-1}(\ol{x}) = \wt{X}$, with $ \ol{x} \in \lbrace \ol{a}, \ol{b}, \ol{c}, \ol{d}, \ol{e}, \ol{f} \rbrace$.
   
   Let $(\wt{C}, \wt{B}), 1 \le \vert \wt{C} \vert \le s, 1 \le \vert \wt{B} \vert \le s$, be a partition of the set of all points incident with some line $L$ of $\mS$. Let $\wt{D}$ consist of the points not incident with $L$ but collinear with a point of $\wt{C}$, let $\wt{A}$ consist of the points not incident with $L$ but collinear with a point of $\wt{B}$, let $\wt{E}$ consist of the points not in $\wt{C} \cup \wt{D}$ but collinear with a point of $\wt{D}$, and let $\wt{F}$ consist of the points not in $\wt{A} \cup \wt{B}$ but collinear with a point of $\wt{A}$. Further, $\Phi^{-1}(\ol{b}\ol{c}) = \lbrace L \rbrace$, $\Phi^{-1}(\ol{c}\ol{d})$ is the set of all lines distinct from $L$ but incident with a point of $\wt{C}$, $\Phi^{-1}(\ol{a}\ol{b})$ is the set of all lines distinct from $L$ but incident with a point of $\wt{B}$, $\Phi^{-1}(\ol{d}\ol{e})$ is the set of all lines distinct from the lines of $\Phi^{-1}(\ol{d}\ol{c})$ but incident with a point of $\wt{D}$, $\Phi^{-1}(\ol{f}\ol{a})$ is the set of all lines distinct from the lines of $\Phi^{-1}(\ol{a}\ol{b})$ but incident with a point of $\wt{A}$, $\Phi^{-1}(\ol{d}\ol{e})$ is the set of all lines distinct from the lines 
   of $\Phi^{-1}(\ol{c}\ol{d})$ but incident with a point of $\widetilde{D}$,
   and $\Phi^{-1}(\ol{f}\ol{e})$ is the set of all lines not in $\Phi^{-1}(\ol{f}\ol{a})$ but incident with a point of $\wt{F}$ (that is, the set of all lines not in $\Phi^{-1}(\ol{e}\ol{d})$ but incident with a point of $\wt{E}).$
   
   \item[{\rm (b)}] The dual of (a).
  
  \end{itemize}
  \end{theorem}
  
  {\it Proof}.\quad As in the proof of the previous theorem, we proceed in a number of steps. We first explain a connection between thin generalized hexagons and projective planes.
  
Let $\mG$ be a thin generalized hexagon of order $(s^\prime, 1), s^\prime \ge 1$. Defines an equivalence relation $\sim$ on the set of lines of $\mG$: \ul{$L \sim M$ if $\widetilde{\bf{d}}$$(L, M)$ = even,} where $\widetilde{\bf{d}}(\cdot,\cdot)$ is the distance in the line graph of $\mG$.  
Let $U, V$ be the equivalence classes. Call the elements of $U$ {\it points}, the elements of $V$ {\it lines}, and $L \in U$ is {\it incident} with $M \in V$ if and only if $\widetilde{\bf{d}}$$(L, M) = 1$. Then this incidence structure is a projective plane $\mP$ of order $s^\prime$ (projective planes of order 1 are allowed). Let $x$ be a point of $\mG$, and let $L \in U$ and $M \in V$ be the lines of $\mG$ which are incident with $x$. Then the point $x$ can be identified with the flag $(L, M)$ of the plane $\mP$.\\
  
\quad(\rm 1) {\bf Let $L$ be any line of $\mS$. Then $\Phi$ maps the set of all points incident with $L$ onto the set of all points incident with $\Phi(L)$}.

{\em Proof of {\em (1)}.}\quad   Let $L$ be any line of $\mS$, let $\Phi(L) = \ol{L}$, let $\widehat{L}$ be the set of all points incident with $L$, and let $\widehat{\ol{L}}$ be the set of all points incident with $\ol{L}$. Assume, by way of contradiction, that $\Phi(\widehat{L}) \not = \widehat{\ol{L}}$. Let $\ol{x} \: \mathbf{I} \: \ol{L}, \ol{x} \not \in \Phi(\widehat{L})$. Let {\bf{d}}$(\ol{x},\ol{y})$ = 4, {\bf{d}}$(\ol{y}, \ol{L})$ = 5 and $y \in \Phi^{-1}(\ol{y})$. If $z \: \mathbf{I} \: L$, {\bf{d}}$(z, y) \le 4$, then {\bf{d}}$(z, y) = 4$, ${\mathbf d}(y,L) = 5$ and $\Phi(z) = \ol{x}$. So $\ol{x} \in \Phi(\widehat{L})$, a contradiction.\\
   
\quad(\rm 2) {\bf Let $x$ be any point of $\mS$, and let $\Phi(x) = \ol{x}$. Then the image of the set of lines incident with $x$, is the set of lines incident with $\ol{x}$}.
   
   {\em Proof of {\rm (2)}}.\quad Dualize the proof of (1).
   
   We will use (1) and (2) without further notice.\\
   
\quad(\rm 3) {\bf We have that $s^\prime = 1$}.
   
{\em Proof of {\rm (3)}}.\quad   Let $A, B$ be the two systems of lines of $\mG$ arising from the aforementioned equivalence relation, let $\Phi^{-1}(A)$ be the union of all inverse images of the elements of $A$, and let $\Phi^{-1}(B)$ be the union of all inverse images of the elements of $B$. If $x$ is a point of $\mS$, then $A_x$ is the set of all lines of $\Phi^{-1}(A)$ incident with $x$ and $B_x$ is the set of all lines of $\Phi^{-1}(B)$ incident with $x$.
   
   Let $c$ be a point of $\mS$ and let $M$ be a line incident with $c$. Let $\Phi(c) = \ol{c}$, let $\ol{d}$ be a point at distance 6 from $\ol{c}$, and let $\Phi(d) = \ol{d}$. Then $\bf{d}$$(c, d) = 6$. The line $\ol{M} = \Phi(M)$ is incident with $\ol{c}$; say $\ol{M} \in B$. Let $\ol{d} \: \mathbf{I} \: \ol{U}, {\bf{d}}(\ol{U}, \ol{M}) = 4$ and $d \: \mathbf{I} \: U, {\bf{d}}(U, M) = 4$. Then $U \in \Phi^{-1}(\ol{U})$ and $\ol{U} \in B$. It easily follows that $\vert B_c \vert = \vert B_d \vert $. Similarly $\vert A_c \vert = \vert A_d \vert $.
   
   From now on, assume, by way of contradiction, that $s^\prime > 1$.
   
   Let $e$ and $f$ be points at distance 6 in $\mS$. If {\bf{d}}$(\ol{e}, \ol{f}) = 6$, with $\Phi(e) = \ol{e}$ and $\Phi(f) = \ol{f}$, then we know that $\vert B_e \vert = \vert B_f \vert$ and $ \vert A_e \vert = \vert A_f \vert $.
   
   Now let {\bf{d}}$(\ol{e}, \ol{f}) = 4$, let  $\ol{e} \: \mathbf{I} \: \ol{E} , \ol{f} \: \mathbf{I} \: \ol{F}$, {\bf{d}}$(\ol{E}, \ol{F}) = 2$ and let {\bf{d}}$(\ol{u}, \ol{E})$ = {\bf{d}}$(\ol{u}, \ol{F}) = 1$. Consider the projective plane $\mP = (A, B)$ of order $s^\prime$ which corresponds with $\mG$;  here $A$ is the point set of $\mP$ and $B$ is the line set of $\mP$. Without loss of generality $\ol{E} \in A, \ol{F} \in B$ and $\ol{u}$ is the flag $(\ol{E}, \ol{F})$. Then $\ol{e}$ is a flag $(\ol{E}, \ol{T})$ and $\ol{f}$ is a flag $(\ol{W}, \ol{F})$. Let $\ol{v} = (\ol{G}, \ol{D})$ be a flag of $\mP$ such that $\ol{G} \not \! \mathbf{I} \: \ol{T}, \ol{E}  \not \! \mathbf{I} \: \ol{D}, \ol{G} \not \! \mathbf{I} \: \ol{F}, \ol{W} \not \! \mathbf{I} \: \ol{D}$ (in $\mP$,  so in $\mG$ we have $\ol{G} \not \sim \ol{T}, \ol{E} \not \sim \ol{D}, \ol{G} \not \sim \ol{F}, \ol{W} \not \sim \ol{D})$; the flag $\ol{v}$ exists as $s^\prime > 1$. Then {\bf{d}}$(\ol{e}, \ol{v})$ = {\bf{d}}$(\ol{f}, \ol{v})$ = 6. So {\bf{d}}$(e, v)$ = {\bf{d}}$(f, v)$ = 6. It follows that $\vert A_e \vert = \vert A_v \vert = \vert A_f \vert $ and $\vert B_e \vert = \vert B_v \vert = \vert B_f \vert. $
   
   Next let {\bf{d}}$(\ol{e}, \ol{f})$ = 2. Consider again the projective plane $\mP = (A, B)$. Without loss of generality $\ol{e}\ol{f} = \ol{N} \in A$. The points $\ol{e}, \ol{f}$ are flags $(\ol{N}, {\ol{N}}_1), (\ol{N}, {\ol{N}}_2)$ of $\mP$. Choose a flag $(\ol{G}, \ol{D}) = \ol{v}$ of $\mP$ such that $\ol{G} \not \! \mathbf{I} \: {\ol{N}}_1, \ol{G} \not \! \mathbf{I} \: {\ol{N}}_2, \ol{N} \not \!\mathbf{I} \: \ol{D} $ (in $\mP$). Then {\bf{d}}$(\ol{e}, \ol{v})$ = {\bf{d}}$(\ol{f}, \ol{v}) = 6$. If $v \in \Phi^{-1}(\ol{v})$, then {\bf{d}}$(e, v)$ = {\bf{d}}$(f, v)$ = 6. So $\vert A_e \vert = \vert A_v \vert = \vert A_f \vert$ and $\vert B_e \vert = \vert B_v \vert = \vert B_f \vert. $

  Now let $\ol{e} = \ol{f}$. Let {\bf{d}}$(\ol{e}, \ol{v})$ = {\bf{d}}$(\ol{f}, \ol{v})$ = 6, with $\ol{v}$ some point of $\mG$, and let $v \in \Phi^{-1}(\ol{v})$. Then {\bf{d}}$(e, v)$ = {\bf{d}}$(f, v)$ = 6 and so $\vert A_e \vert = \vert A_v \vert = \vert A_f \vert$ and $\vert B_e \vert = \vert B_v \vert = \vert B_f \vert. $
  
  Hence for {\bf{d}}$(e, f)$ = 6 we always have $\vert B_e \vert = \vert B_f \vert$ and $\vert A_e \vert = \vert A_f \vert$ .
  
  Next, consider points $e$ and $f$ of $\mS$ at distance 4 or 2. Clearly we can find a point $v$ of $\mS$ with {\bf{d}}$(e, v)$ = {\bf{d}}$(f, v)$ = 6. Then $\vert A_e \vert = \vert A_v \vert = \vert A_f \vert$ and $\vert B_e \vert = \vert B_v \vert = \vert B_f \vert$.
  
  We conclude that for any two points $e, f$ of $\mS$ we have $\vert A_e \vert = \vert A_f \vert$ and $\vert B_e \vert = \vert B_f \vert$.
  
  Let $r$ be a point of $\mS$. At least one of $\vert A_r \vert, \vert B_r \vert$ is at least two; say $\vert A_r \vert \ge 2$. Assume that the lines $L_1, L_2, L_1 \not = L_2$, of $\mS$ belong to $A_r$. Then $\Phi(L_1) = \Phi(L_2) = \ol{L} \in A$ and $\ol{r} \: \mathbf{I} \: \ol{L}$, with $\ol{r} = \Phi(r)$. Let $\ol{m} \: \mathbf{I} \: \ol{L}, \ol{m} \not = \ol{r}$, and $m \: \mathbf{I} \: L_2$, with $m \in \Phi^{-1}(\ol{m})$. Let $m \: \mathbf{I} \: N$, with $N \in B_m$. Choose $u_1 \: \mathbf{I} \: L_1, u_1 \not = r$ and $u_2 \: \mathbf{I} \: N, u_2 \not = m$, with $\Phi(u_1) \not = \ol{m} \not = \Phi(u_2)$. Then {\bf{d}}$({u_1},{u_2}) = 6$. Let $u_1 \: \mathbf{I} \: R_1$, with $R_1 \in B_{u_1}$, and {\bf{d}}$(R_1, R_2) = 4$ with $u_2 \: \mathbf{I} \: R_2$. Hence $R_1$ determines $R_2$. Then $\Phi(R_2) = {\ol{R}}_2 = \ol{N} \in B$. So $R_2 \in B_{u_2}$. It follows that $\vert B_{u_1} \vert \leq \vert B_{u_2} \setminus \lbrace N \rbrace \vert$. So $\vert B_{u_1} \vert < \vert B_{u_2} \vert$, a contradiction.
  
  Consequently $s^\prime = 1$.\\
  
 \quad{\rm (4)} {\bf Let $s^\prime = 1$. Then not for any given points $\ol{x}$ and $\ol{y}$ in $\mG$, with {\bf{d}}$(\ol{x}, \ol{y}) = 4$, there exist points $x \in \Phi^{-1}(\ol{x})$ and $y \in \Phi^{-1}(\ol{y})$ with {\bf{d}}$(x, y) = 6$}.
  
Let {\bf{d}}$(\ol{u}, \ol{v})$ = 6, with $\ol{u} = \Phi(u)$ and $\ol{v} = \Phi(v)$, be points of $\mG$. Then by (3) we have $\vert B_u \vert = \vert B_v \vert$ and $\vert A_u \vert = \vert A_v \vert$.\\
  
 {\em Proof of {\rm (4)}}.\quad Assume, by way of contradiction, that for any given points $\ol{x}$ and $\ol{y}$ in $\mG$, with {\bf{d}}$(\ol{x}, \ol{y})$ = 4, there exist points $x \in \Phi^{-1}(\ol{x})$ and $y \in \Phi^{-1}(\ol{y})$ with {\bf{d}}$(x, y)$ = 6.
  
  Let $\ol{a}, \ol{b}, \ol{c}, \ol{d}, \ol{e}, \ol{f}$ be the points of $\mG$, with $\ol{a} \sim \ol{b} \sim \ol{c} \sim \ol{d} \sim \ol{e} \sim \ol{f} \sim \ol{a}$. Let $\ol{a} = \Phi(a)$. Choose $f, e, d, N $ such that $\ol{f} = \Phi(f), \ol{e} = \Phi(e) , \ol{d} = \Phi(d), \ol{d}\ol{c} = \Phi(N)$ and $d \: \mathbf{I} \: N$, with $a \sim f \sim e \sim d$. Let $c \: \mathbf{I} \: N$ with {\bf{d}}$(a, c)$ = 4, and let {\bf{d}}$(a, b)$ = {\bf{d}}$(b, c)$ = 2. Then $\ol{b} = \Phi(b)$ and $\ol{c} = \Phi(c)$. Let $A = \lbrace \ol{a}\ol{f}, \ol{b}\ol{c}, \ol{e}\ol{d} \rbrace$ and $B = \lbrace \ol{a}\ol{b}, \ol{c}\ol{d}, \ol{e}\ol{f} \rbrace$. We have $\vert B_a \vert = \vert B_d \vert, \vert B_b \vert = \vert B_e \vert, \vert B_f \vert = \vert B_c \vert$ and $\vert A_a \vert = \vert A_d \vert, \vert A_b \vert = \vert A_e \vert, \vert A_f \vert = \vert A_c \vert$.
  
  Let $\Phi(a^\prime) = \ol{a}$. As {\bf{d}}$(\ol{a}, \ol{d})$ = 6, we have {\bf{d}}$(a^\prime, d)$ = 6. Hence $\vert B_a \vert = \vert B_d \vert = \vert B_{a^\prime} \vert$, and similarly $\vert A_a \vert = \vert A_{a^\prime} \vert$.
  
  By assumption there is a point $c^\prime$, with $\Phi(c^\prime) = \ol{c}$ and {\bf{d}}$(a, c^\prime)$ = 6. Let $c^\prime \sim n \sim m \sim a$ and choose the point $m$ in such a way that $\Phi(ma) = \ol{f}\ol{a}$. Then $\Phi(m) = \ol{a}$ and $\Phi(n) = \ol{b}$. So $\Phi(nc^\prime) = \ol{b}\ol{c} \in A$. Hence $\vert A_{c^\prime} \vert \ge \vert A_a \vert$, so $\vert A_c \vert \ge \vert A_a \vert$. Also $\vert A_c \vert = \vert A_f \vert, \vert A_a \vert = \vert A_d \vert$, and so $\vert A_c \vert \ge \vert A_d \vert$.
  
  Similarly $\vert A_d \vert \ge \vert A_c \vert$, and so $\vert A_c \vert = \vert A_d \vert$.
  
  Consequently $\vert A_a \vert = \vert A_b \vert = \vert A_c \vert = \vert A_d \vert = \vert A_e \vert = \vert A_f \vert$ and $\vert B_a \vert = \vert B_b \vert = \vert B_c \vert = \vert B_d \vert = \vert B_e \vert = \vert B_f \vert$. Hence $\vert A_u \vert $ and $\vert B_u \vert$ are independent of the choice of the point $u$ of $\mS$.
  
  Such as in the last part of (3) we now find a contradiction.\\
  
\quad{\rm (5)} {\bf For $s^\prime = 1$ two (mutually dual) cases occur}.
  
{\em Proof of {\rm (5)}}.\quad  By (4) we know that not for any given points $\ol{x}$ and $\ol{y}$ in $\mG$ with {\bf{d}}$(\ol{x}, \ol{y})$ = 4, there exist points $x \in \Phi^{-1}(\ol{x})$ and $y \in \Phi^{-1}(\ol{y})$ with {\bf{d}}$(x, y)$ = 6.
  
  Let $\ol{a}, \ol{b}, \ol{c}, \ol{d}, \ol{e}, \ol{f}$ be the points of $\mG$, with $\ol{a} \sim \ol{b} \sim \ol{c} \sim \ol{d} \sim \ol{e} \sim \ol{f} \sim \ol{a}$, let $A = \lbrace \ol{b}\ol{c}, \ol{d}\ol{e}, \ol{f}\ol{a} \rbrace$ and $B = \lbrace \ol{a}\ol{b}, \ol{c}\ol{d}, \ol{e}\ol{f} \rbrace$, and let $\Phi^{-1}(\ol{x}) = \wt{X}$ for $\ol{x} = \lbrace \ol{a}, \ol{b}, \ol{c}, \ol{d}, \ol{e}, \ol{f} \rbrace$.
  
  Without loss of generality, we assume that there do not exist points $a^\prime \in \wt{A}$ and $c^\prime \in \wt{C}$ for which {\bf{d}}$(a^\prime, c^\prime)$ = 6 (and so {\bf{d}}$(a^\prime, c^\prime)$ = 4).\\

\quad {\bf CASE (a)} {\it Assume that $b_1 \sim b_2, b_1 \not = b_2$, with $b_1, b_2 \in \wt{B}$}.

We distinguish two cases. 

\begin{itemize}
\item[(a.1)] {\it We have that $\Phi(b_1b_2) = \ol{b}\ol{c}  \in A$}.
  
  Let $M \in \Phi^{-1}(B)$ be a line incident with $b_1$ which is incident with some point $a_1 \in \wt{A}$. Assume, by way of contradiction, that the line $M^\prime \not = b_1b_2$ is incident with $b_2$ and belongs to $\Phi^{-1}(A)$. Then $M^\prime$ is incident with some point $c_1$ of $\wt{C}$. As {\bf{d}}$(a_1, c_1)$ = 6, we have a contradiction. Hence $b_2$, and similarly $b_1$, is incident with $t$ lines of $\Phi^{-1}(B)$.
  
  Assume, by way of contradiction, that there exists a point $b^\prime \in \wt{B}$, with $b^\prime \sim b_1$, and $b^\prime \not \! \mathbf{I} \; b_1b_2$. Then, interchanging $b^\prime$ and $b_2$, and $A$ and $B$, we see that $b_1$ and $b^\prime$ are incident with $t$ lines of $\Phi^{-1}(A)$. Hence $t = 1$, a contradiction. Hence each line $M \in \Phi^{-1}(B)$ incident with $b_1$ is incident with exactly $s$ points of $\wt{A}$. Similarly, each line of $\Phi^{-1}(B)$ incident with some $b_i \in \wt{B}, b_i \: \mathbf{I} \: b_1b_2$, is incident with $s$ points of $\wt{A}$.
  
  Let $b_1, b_2, \ldots, b_u$ be the points of $\wt{B}$ incident with $b_1b_2$. Then all points not incident with $b_1b_2$ but collinear with one of $b_1, b_2, \ldots, b_u$ are elements of $\wt{A}$. 
  
  Assume that $u < s$, so $b_1b_2$ is incident with at least  two points $c_1, c_2$ of $\wt{C}$. If $a^\prime \in \wt{A}$ and $a^\prime$ is collinear with no point of $\lbrace b_1, b_2, \ldots, b_u \rbrace$, then {\bf{d}}$(a^\prime, c_1)$ = {\bf{d}}$(a^\prime, c_2)$ = 4, yielding a contradiction. Hence for $u < s$ the set $\wt{A}$ consists of all points not incident with $b_1b_2$ but collinear with one of $b_1, b_2, \ldots, b_u$.
  
Let $u \leq s$ again.  Assume, by way of contradiction, that $c^\prime \in \wt{C}, c^\prime \not \! \mathbf{I} \: b_1b_2$. Let $b_1 \sim a^{\prime\prime} \sim a^{\prime\prime\prime} \sim b_1$, with $a'', a''' \in \widetilde{A}$, 
$a^{\prime\prime} \not = a^{\prime\prime\prime}, a^{\prime\prime} \not = b_1, a^{\prime\prime\prime} \not = b_1$.  
  Then {\bf{d}}$(c^\prime, a^{\prime\prime})$ = {\bf{d}}$(c^\prime, a^{\prime\prime\prime})$ = 4 implies that $c^\prime \sim b_1$. Similarly $c^\prime \sim b_2$, a contradiction. Consequently all points of $\wt{C}$ are incident with the line $b_1b_2$.
  
  Assume that there is a point $b^{\prime\prime} \not \! \mathbf{I} \: b_1b_2, b^{\prime\prime} \in \wt{B}$. A line $V$ of $\Phi^{-1}(A)$ incident with $b^{\prime\prime}$ is incident with a point $c^{\prime\prime}$ of $\wt{C}$ ($c^{\prime\prime} \: \mathbf{I} \: b_1b_2$); the other $s$ points incident with $V$ belong to $\wt{B}$.\\
  
\quad(A) {\it Assume that $\Phi^{-1}(\ol{b}\ol{c}) = \lbrace b_1b_2 \rbrace$}.
  
  Assume, by way of contradiction, that $b^{\prime\prime} \in \wt{B}$ with $b^{\prime\prime} \not \! \mathbf{I} \: b_1b_2$. Then $b^{\prime\prime}$ is incident with a line $W$ of $\Phi^{-1}(A) \setminus \lbrace b_1b_2 \rbrace$. Hence $W \in \Phi^{-1}(\ol{b}\ol{c}), W \not = b_1b_2$, a contradiction. Hence each point of $\wt{B}$ is incident with the line $b_1b_2$. In this case $\wt{B}, \wt{C}$ determine a partition of the set of points incident with $b_1b_2, 1 \le \vert \wt{B} \vert \le s, 1 \le \vert \wt{C} \vert \le s$, the set $\wt{A}$ consists of all points not incident with $b_1b_2$ but collinear with a point of $\wt{B}$, and $\wt{D}$ consists of all points not incident with $b_1b_2$ but collinear with a point of $\wt{C}$. Now it is easy to see that we have Case (a) in the statement of the theorem.\\
  
 \quad(B) {\it Assume that $\vert \Phi^{-1}(\ol{b}\ol{c}) \vert > 1$}.
  
  If $\vert \wt{C} \vert > 1$, then we know that $\wt{A}$ consists of all points not incident with $b_1b_2$ but collinear with one of the points of $\wt{B}$ incident with $b_1b_2$. Let $U \in \Phi^{-1}(\ol{b}\ol{c}), U \not = b_1b_2$, let $b^{\prime\prime} \: \mathbf{I} \: U$, with $b^{\prime\prime} \in \wt{B}$, let $U^\prime$ be a line of $\Phi^{-1}(B)$ incident with $b^{\prime\prime}$, and let $a^{\prime\prime} \in \wt{A}$ be a point incident with $U^\prime$. Then $a^{\prime\prime}$ must be collinear with a point of $\wt{B}$ incident with $b_1b_2$, clearly a contradiction.\\
  
  Hence $\wt{C} = \lbrace c_1 \rbrace$. In such a case $\wt{B}$ consists of all points distinct from $c_1$ on $r$ lines incident with $c_1$, $1 \le r \le t$, and $\wt{D}$ consists of all points distinct from $c_1$ incident with $t + 1 - r$ lines incident with $c_1$. Now it is easy to see that we have Case (b) (that is the dual of Case (a)) in the statement of the theorem.\\

\item[(a.2)] {\it We have that $\Phi(b_1b_2) \in B$}.\\

  This case is similar to Case (a.1).
\end{itemize}

\quad{\bf CASE (b)} {\it No two points of $\wt{B}$ are collinear}.

  Assume, by way of contradiction, that $b_1, b_2 \in \wt{B}, b_1 \not = b_2$. Let $M$ be a line of $\Phi^{-1}(B)$ incident with $b_2$, and let $N$ be a line of $\Phi^{-1}(A)$ incident with $b_1$. Let $a_1, a_2$ be distinct points of $\wt{A}$ incident with $M$ ($M$ is incident with $s$ points of $\wt{A}$), and let $c_1$ be a point of $\wt{C}$ incident with $N$ ($N$ is incident with $s$ points of $\wt{C}$). As {\bf{d}}$(a_1, c_1)$ = {\bf{d}}$(a_2, c_1)$ = 4, we have $c_1 \sim b_2$. If $c_2 \in \wt{C} \setminus \lbrace c_1 \rbrace$ is incident with $N$, then similarly $c_2 \sim b_2$, yielding a contradiction. Hence $\wt{B} = \lbrace b_1 \rbrace$,
  
  Then $\wt{A}$ consists of the points distinct from $b_1$ incident with $r^\prime$ lines incident with $b_1, 1 \le r^\prime \le t$, and $\wt{C}$ consists of the points distinct from $b_1$ incident with $t + 1 - r^\prime$ lines incident with $b_1$. 
  
  Now it is easy to see that we have Case (b) (that is the dual of Case (a)) in the statement of the theorem.\eop \\

  \section{Locally finitely chained and generated generalized polygons}
  \label{chain}

Recall that a group $G$ is called {\em locally finite} if every finite subset of $G$ generates a finite subgroup. 

\subsection*{EXAMPLES} 

Finite groups and infinite direct sums of finite groups are examples. Also, every subgroup of a locally finite group is evidently locally finite. A {\em Dedekind group} is a group in which each subgroup is normal (example: abelian groups). Call a group {\em Hamiltionian} if it is Dedekind and nonabelian. Then every Hamiltionian group is locally finite. Many other interesting examples exist. \\ 

We now introduce the analogons for generalized polygons. \\

Let $S$ be any point set of a generalized $n$-gon $\Gamma$; then $\langle S \rangle$ by definition is the intersection of all (thin and thick) sub $n$-gons that contain $S$. We call $\langle S \rangle $ the subgeometry {\em generated by $S$}. Note that this not necessarily is a (thin or thick) generalized $n$-gon itself. But generically it is. 
 
Call a thick generalized $n$-gon {\em locally finitely generated} if the following property holds: 
\begin{quote}
for every finite point subset $S$ which generates a (possibly thin) sub $n$-gon $\langle S \rangle$, we have that $\langle S \rangle$ is finite.  
\end{quote}
Finite generalized $n$-gons are trivially locally finitely generated. \\

Call a generalized $n$-gon $\Gamma$ {\em locally finitely chained} if there exists a chain of finite point subsets
\[  S_0 \subseteq S_1 \subseteq \ldots \subseteq S_i \subseteq \ldots                    \]       
indexed over the positive integers, such that each $S_j$ generates a finite sub $n$-gon (possibly thin), and such that 
\[  \bigcup_{i \geq 0}\Big\langle S_i \Big\rangle\ =\ \Gamma.  \]
If $\Gamma$ is not finite, it is easy to see that we may assume that the chain 
\[  \langle S_0 \rangle \subseteq \langle S_1\rangle \subseteq \ldots \subseteq \langle S_i \rangle \subseteq \ldots                    \]     
is strict.\\

A very simple key observation is the following lemma, in which we suppose that locally finitely generated and locally finitely 
chained $n$-gons are not finite (to avoid trivialities). 

\begin{lemma}
\label{lem1}
\begin{itemize}
\item[{\rm (a)}]
Every thick locally finitely generated generalized $n$-gon contains thick locally finitely chained sub $n$-gons. 
\item[{\rm (b)}]
Every thick finitely chained generalized $n$-gon is locally finitely generated.
\item[{\rm (c)}]
A thick generalized $n$-gon is locally finitely generated if and only if its point-line dual is.
\item[{\rm (d)}]
A thick generalized $n$-gon is locally finitely chained if and only if its point-line dual is.
\end{itemize}
\end{lemma}

{\em Proof}.\quad
(a) Consider a chain of finite point subsets 
\[  S_0 \subseteq S_1 \subseteq \ldots \subseteq S_i \subseteq \ldots                    \]     
such that each $S_j$ generates a generalized $n$-gon and not each $\langle S_j \rangle$ is thin, and such that if $i \ne j$ then 
$\langle S_i \rangle \ne \langle S_j \rangle$. (Note that such chains exist!) 
 Then $\cup_{i \geq 0}\Big\langle S_i \Big\rangle$ is thick, not finite  and locally finitely chained. \\

(b) Let $\Gamma$ be thick and locally finitely chained. Let
\[  S_0 \subseteq S_1 \subseteq \ldots \subseteq S_i \subseteq \ldots                    \]    
be a chain of finite point subsets such that each $\langle S_i \rangle$ is a finite generalized $n$-gon and so that 
$\cup_{i \geq 0}\Big\langle S_i \Big\rangle = \Gamma$. Then not each $\langle S_j \rangle$ is thin. 
We assume that $\Gamma$ is not finite. 
Consider any finite point subset $S$ of $\Gamma$. 
Then there is some index $u$ such that $S \subseteq \langle S_u \rangle$. It follows that $\langle S \rangle$ is finite. \\

(c) and (d) are obvious.
\eop \\

Note that $\cup_{i \geq 0}\Big\langle S_i \Big\rangle$ is countably infinite or finite.\\

\begin{observation}
\begin{itemize}
\item[{\rm (a)}]
Let $\Gamma$ be a thick locally finitely chained generalized $n$-gon. Then $n \in \{ 3, 4, 6, 8 \}$.
\item[{\rm (b)}]
Let $\Gamma$ be a thick locally finitely generated generalized $n$-gon. Then $n \in \{ 3, 4, 6, 8 \}$.
\end{itemize} 
\end{observation}

{\em Proof}.\quad
(a) Suppose 
\[  S_0 \subseteq S_1 \subseteq \ldots \subseteq S_i \subseteq \ldots                    \]     
is a chain of finite point subsets for which $\langle S_j \rangle$ is a finite generalized $n$-gon for each $j$, and such that $\cup_{i \geq 0}\Big\langle S_i \Big\rangle = \Gamma$. Then for some $n \in \mathbb{N}$, we have that $\langle S_m \rangle$ is thick as soon as $m \geq n$. For such an $m$, $\langle S_m \rangle$ is a finite thick generalized $n$-gon, so by Feit and Higman \cite[Theorem 1.7.1]{POL}, $n \in \{ 3, 4, 6, 8   \}$.\\

(b) Follows from part (a) and Lemma \ref{lem1}. \eop \\

Before proceeding, we introduce the following two properties for a point-line incidence structure $\Gamma = (\mP,\mL,\I)$:
\begin{itemize}
\item[(Ia)]
every two elements of $\mP \cup \mL$ can be joined  by at most one path of length $< n$; 
\item[(Ib)]
every two elements of $\mP \cup \mL$ can be joined by at least one path of length $\leq n$. 
\end{itemize}

Call a point-line incidence geometry $\Gamma = (\mP,\mL,\I)$ {\em firm} if each element is incident with at least two different elements. 
Then by Van Maldeghem \cite{POL2}, $\Gamma$ is a weak generalized $n$-gon if and only if it is firm, and both (Ia) and (Ib) are satisfied. 

\begin{theorem}
Let $\Gamma$ be a thick locally finitely chained generalized $n$-gon with $n \in \{ 4, 6 \}$, and let $\gamma: \Gamma \mapsto \Delta$ be an epimorphism, with $\Delta$ a thin $n$-gon of order $(s,1)$. Then $s = 1$ and the conclusions of Theorem \ref{JATGQ} and Theorem \ref{JATGH} hold.
\end{theorem}

{\em Proof}.\quad
{\bf CASE $n = 4$}\\
We use the notation introduced in this section. Assume, by way of contradiction, that 
 $\Delta'$ is a finite subgrid of $\Delta$ of order $(s',1)$ with $s' > 1$. 
As $\Gamma = \cup_{i \geq 0}\Big\langle S_i \Big\rangle$, we know that there is some positive integer $m$ for which $\Delta' \subseteq 
\gamma(\Big\langle S_m \Big\rangle)$. Since $\Big\langle S_m \Big\rangle$ is a generalized quadrangle, it follows that $\gamma(\Big\langle S_m \Big\rangle) =: \Delta''$ is a finite $(u \times v)$-grid, with $u, v \geq s'$. Also, 
\[ \gamma:\ \Big\langle S_m \Big\rangle\ \mapsto\ \Delta''              \]
is an epimorphsim. By the proof of Theorem \ref{JATGQ} (in the proof of Theorem \ref{JATGQ} we have $u = v$ but it generalizes trivially), the statement follows.\\

{\bf CASE $n = 6$}\\
We follow essentially the same procedure as in the CASE $n = 4$ (and the same notation). Assume, by way of contradiction, that $s \ne 1$. 
Now define $m$ such that $\gamma(\Big\langle S_m \Big\rangle)$ strictly contains an ordinary $n$-gon $\Delta'$ as subgeometry, and such that $\Big\langle S_m \Big\rangle$ is thick. Since  $\gamma(\Big\langle S_m \Big\rangle) =: \Delta''$ is contained in $\Delta$, we have that (Ia) is automatically satisfied in $\Delta''$. And (Ib) can be verified in $\Delta''$ by applying $\gamma^{-1}$. Since $\Delta''$ contains an ordinary $n$-gon, it now follows that it is firm, so it is a weak generalized $n$-gon. Also, 
\[ \gamma:\ \Big\langle S_m \Big\rangle\ \mapsto\ \Delta''              \]
is an epimorphsim. By \cite[Theorem 3.1]{POL2} (see also \cite[Theorem 1.6.2]{POL}), $\Delta''$ is either the point-line dual of the double of a projective plane, or the point-line dual of the double of a degenerate plane $\mP$ which consists of a nonincident point-line pair $(a,A)$, and at least two lines incident with $a$ which contain one other point, on $A$. The point-line dual of the double of $\mP$ has at most two thick lines, so we can take $m$ large enough so that $\gamma(\langle S_m \rangle)$ contains more than two thick lines.
By the proof of Theorem \ref{JATGH}, the statement follows (as $\Delta''$ is the point-line dual of the double of a projective plane).\\
 \eop \\
 

\begin{corollary}
Let $\Gamma$ be a thick locally finitely generated generalized $n$-gon with $n \in \{ 4, 6 \}$, and let $\gamma: \Gamma \mapsto \Delta$ be an epimorphism, with $\Delta$ a thin $n$-gon of order $(s,1)$. Then $s = 1$ and the conclusions of Theorem \ref{JATGQ}  and Theorem \ref{JATGH} hold. \eop \\
\end{corollary}

\subsection{Examples}

Although at first sight the property of being locally finitely chained / generated appears to be rather strong, there are surprisingly many interesting examples of such polygons with much structure. \\

Suppose $\F_q$ is a finite field, and let $\overline{\F_q}$ be an algebraic closure of $\F_q$. Let $q = p^m$ for the prime $p$. Then for each positive integer $n$, $\overline{\F_q}$ has a unique subfield isomorphic to $\F_{p^n}$ (we denote it also as $\F_{p^n}$).  Then 
\[ \bigcup_{i \geq 1}\F_{p^i}\ =\ \overline{\F_q}.    \]

\subsection*{Classical examples of quadrangles}

Now consider the quadratic equation
\begin{equation}
\label{eq1}   X_0^2 + X_1X_2 + X_3X_4 = 0.            \end{equation}

Fix a prime $p$. For each natural power $q$ of $p$, the $\F_q$-rational points and $\F_q$-rational lines of the variety defined by eq. (\ref{eq1}) define a generalized quadrangle $\mQ(4,q)$ of order $(q,q)$ (we work in a projective space $\bP^4(q)$ with homogeneous coordinates $(X_0: X_1: X_2: X_3: X_4)$.) Obviously, if we consider the equation over $\overline{\F_p} = \overline{\F_q}$, we obtain a generalized quadrangle which lives in $\bP^4(\overline{\F_p})$, denoted $\mQ(4,\overline{\F_p})$, and for which
\[    \mQ(4,\overline{\F_p})  = \bigcup_{i \geq 0}   \mQ(4,{\F_{p^i}}).      \]

The quadrangle $\mQ(4,\overline{\F_p})$ is locally finitely chained (and as it is countably infinite, it is also locally finitely generated).

\subsection*{Nonclassical examples of quadrangles}

\subsubsection{Direct limit of ovals}

The example in this section is taken from \cite{KTGroring}.
We first recall a result of Segre. 

\begin{theorem}[Segre]
\label{segre}
Let $i$ and $h$ be positive integers such that $\mathrm{gcd}(i,h) = 1$, and let $q = 2^h$. Then the set of points in $\bP^2(q)$ given 
by homogeneous coordinates
\begin{equation}
\{ (1 : t : t^{2^i})\ \vert \ t \in \F_q \}\ \cup\ \{ (0 : 0 : 1)   \}
\end{equation}
is an oval $\mO(i,h)$ with nucleus $\eta = (0 : 1 : 0)$. 
\end{theorem}

If $\mathrm{gcd}(i,h) \ne 1$, then the set $\{ (1 : t : t^{2^i})\ \vert \ t \in \F_q \} \cup \{ (0 : 1 : 0), (0 : 0 : 1)\}$ does not form a hyperoval; see Hirschfeld \cite[\S 8.4, Corollary 8.26]{Hirsch}. Each of the ovals constructed from Theorem \ref{segre} has as nucleus $\eta = (0 : 1 : 0)$ (in planes represented over different fields, that is to say). 

Now define a direct system of fields, as follows. First of all, fix one positive integer $i > 1$ (for the sake of convenience, this could be some prime). Now let $\mS(i)$ 
be the set of all positive integers $r$ such that $\mathrm{gcd}(r,i) = 1$. We define a directed set $\Big(\mS(i),\preceq\Big)$ as follows:  $a \preceq  b$ if $a$ divides $b$
(which happens if and only if $\F_{2^a}$ is a subfield of $\F_{2^b}$). We then consider the family $\{ \mO(i,h)\ \vert \ h \in \mS(i) \}$. If $u \preceq v$, 
then define the natural embeddings
\begin{equation}
\begin{cases}
&\iota^\mO_{uv} := \mO(i,u) \ \hookrightarrow\ \mO(i,v), \\
&\iota^\F_{uv}  := \F_{2^u} \ \hookrightarrow\ \F_{2^v}, \\
&\iota^\bP_{uv} := \bP^2(2^u) \ \hookrightarrow\ \bP^2(2^v).
\end{cases}
\end{equation}

Obviously $\Big( \mO(i,u),\iota^\mO_{uv} \Big)$ is a direct system over $\mS(i)$,   as well as $\Big( \F_{2^u},\iota^\F_{uv} \Big)$ and  $\Big( \bP^2(2^u),\iota^\bP_{uv} \Big)$. 

\begin{theorem}[K. Thas \cite{KTGroring}]
We have that $\varinjlim\mO(i,u)$ is an oval of $\bP^2(\ell)$, with nucleus $(0: 1 :0)$; here, $\ell = \varinjlim\F_{2^u}$ and $\bP^2(\ell) = \varinjlim\bP^2(2^u)$. 
\end{theorem}

\medskip
\subsubsection{The quadrangles}

Passing to the corresponding GQs (by the standard $\hT_2(\mO)$-construction \cite[chapter 8]{PT2}), we obtain a direct system of GQs $\Gamma(\mO(i,u),\eta)$ over the directed set $\mS(i)$, and 
\begin{equation}
\varinjlim \Gamma(\mO(i,u),\eta) = \Gamma(\widetilde{\mO},\eta), 
\end{equation}
with $\widetilde{\mO} = \varinjlim \mO(i,u)$ an oval of $\mathbb{P}^2(\ell)$ and $\ell = \varinjlim \F_{2^u}$.\\

Note that 
\[ \Gamma(\widetilde{\mO},\eta) = \bigcup_{u \in \mS(i)}\Gamma(\mO(i,u),\eta)\ \mbox{and}\ \ell = \bigcup_{\mathrm{gcd}(i,v) = 1}\F_{2^v}.                            \]
(Here, each $\F_{2^v}$ is the unique subfield of $\ell$ isomorphic to $\F_{2^v}$.)\\


The quadrangle $\Gamma(\widetilde{\mO},\eta)$ is locally finitely chained (and as it is countably infinite, it is also locally finitely generated).


\subsection*{Examples of hexagons}

Let $k = \F_q$ be a finite field, with $q$ a power of the prime $p$. On the parabolic quadric $\mQ(6,k)$ in $\mathbb{P}^6(k)$ with defining equation
\begin{equation}
X_0X_4 + X_1X_5 + X_2X_6 = X_3^2,
\end{equation}  
we define the generalized hexagon $\mH(k)$ (we work with homogeneous coordinates $(X_0 : X_1 : X_2 : X_3 : X_4 : X_5 : X_6)$ in $\mathbb{P}^6(k)$). 
The lines of $\mH(k)$ are defined through their Grassmann coordinates, which have to satisfy six linear equations with coefficients in $\F_p$ | see \cite[2.4.13]{POL}. A similar construction as in the example $\mQ(4,\overline{\F_p})$ in the beginning of this section, shows that the hexagon 
\[    \mH(\overline{\F_p})  = \bigcup_{i \geq 0}   \mH({\F_{p^i}})      \]
is locally finitely chained (and as it is countably infinite, it is also locally finitely generated).



\medskip
\section{Counter examples in the infinite case}
\label{counter}

In this section we show that Theorem \ref{JATGQ} does not hold for infinite generalized quadrangles, by constructing epimorphisms from 
thick infinite generalized quadrangles to grids with parameters $(s',1)$, $s' > 1$. \\

Suppose that $A$ is a point-line configuration without ordinary triangles as subgeometry, and let $\epsilon: A \mapsto \Gamma$ be an 
epimorphism (of point-line incidence geometries), where $\Gamma$ is a generalized quadrangle of order $(s',1)$, where $s' > 1$ is arbitrary but finite or countable. We also suppose that $A$ has a finite or countable number of points and/or lines.

If $A$  already is a generalized quadrangle, there is nothing to prove, so suppose it is not. 

We will freely construct a generalized quadrangle $\overline{A}$ over $A$ and an epimorphism $\overline{\epsilon}: \overline{A} \mapsto \Gamma$. In each next step, we look at the nonincident point-line pairs $(u,U)$ for which there is no line on $u$ meeting $U$, and add a new incident point-line pair $(v,V)$ for which $u \I V \I v \I U$ (so both $v$ and $V$ are added in this step). 

Let $A = A_1$.

Suppose in a previous step, $A_n$ and $\epsilon_n: A_n \mapsto \Gamma$ were constructed (with $\epsilon_n$ an epimorphism). 

First suppose the point $x$ was added in the previous step, and let $\epsilon_n(x) = \widetilde{x} \in \Gamma$. Let $X$ be a line 
of $A_n$, not incident with $x$ and  such that there is no line yet incident with $x$ which hits $X$. (Note that there are precisely two lines incident with $x$ in $A_n$.)
Then add the new line $Y \I x$ and the new point $y$ for which $Y \I y \I X$. 
If $\epsilon_n(X)$ is not incident with $\widetilde{x}$ in $\Gamma$, then $\epsilon_{n + 1}(y)$ and $\epsilon_{n + 1}(Y)$ are uniquely determined.
If $\epsilon_n(X) \I \widetilde{x}$, then let $\epsilon_{n + 1}(Y) = \epsilon_{n + 1}(X)$ and $\epsilon_{n + 1}(y) = \widetilde{x}$ (and note that other choices are possible). 

Proceed in a similar way for lines $X$ which were added in $A_n$. 

Now let $\overline{A} = \cup_{n \geq 1}A_n$ and $\overline{\epsilon} = \cup_{n \geq 1}\epsilon_n$ (where we identify a function with its graph). Then 
\[  \overline{\epsilon}:\ \overline{A}\ \mapsto\ \Gamma      \]
is an epimorphism. \\

Note that it is easy to construct legitimate  begin-configurations $(A,\epsilon)$:  for instance, let $A$ be $\Gamma$ with one extra isolated point $w$, and let $\epsilon$ be defined as the identity on $\Gamma$, and give $w$ an arbitrary point image. \\

Obviously, this kind of construction can be generalized to other gonalities, and with more general targets. In Gramlich and Van Maldeghem \cite{POL2}, it is proved that starting from any (thick) generalized $n$-gon $\Gamma$ ($n \geq 2$), one can construct a free generalized $n$-gon $\overline{A}$ and an epimorphism 
\[  \overline{\epsilon}:\ \overline{A}\ \mapsto\ \Gamma.      \]

The begin configuration in \cite{POL2} is different, and the target is thick, 
but of course the essence is the same.

\medskip
\section{Remark: epimorphisms involving locally finite quadrangles}
\label{lf}

Call a thick generalized $n$-gon of order $(s,t)$ {\em locally finite} if exactly one of $s, t$ is finite. 

\begin{theorem}
Let $\gamma: \Gamma \mapsto \Delta$ be an epimorphism from the thick locally finite generalized quadrangle $\Gamma$ of order $(s,t)$ with $t$ finite, to the grid $\Delta$ of order $(s',1)$. Then $s' = 1$ and we have the conclusion of Theorem \ref{JATGQ}.  
\end{theorem}

{\em Proof}.\quad
The part of the proof of Theorem \ref{JATGQ} that shows that $s' = 1$ only uses the fact that the number of lines incident with a point is finite. \eop \\ 

Similarly, we have the following result.

\begin{theorem}
Let $\gamma: \Gamma \mapsto \Delta$ be an epimorphism from the thick locally finite generalized hexagon $\Gamma$ of order $(s,t)$ with $t$ finite, to the generalized hexagon $\Delta$ of order $(s',1)$. Then $s' = 1$ and we have the conclusion of Theorem \ref{JATGH}.  
\eop
\end{theorem}



\medskip
\section{Epimorphisms between thin polygons}
\label{thin}

In this section, we have a brief look at epimorphisms with target a thin generalized $n$-gon (with an order), and where the source is also thin. In this entire section, we only work with finite polygons. 


\subsection{Doubling epimorphisms}

Now let $\gamma: \Gamma^{\Delta} \mapsto {\Gamma'}^\Delta$ be an epimorphism, where $\Gamma = (\mP, \mL, \I)$ and $\Gamma' = (\mP', \mL', \I')$ are generalized $n$-gons of order $(s,s)$, respectively $(s',s')$ ($\mP \cup \mL$ and $\mP' \cup \mL'$ are the lines of $\Gamma^\Delta$ and ${\Gamma'}^\Delta$). Suppose $\gamma$ sends some point $v$ in $\mP$ to a point $v'$ in $\mP'$. Let $w \sim v \ne w$ be a point of $\mP$, and let $U = vw$. In $\Gamma$, to $v$ and $w$ correspond lines $V$ and $W$ of $\Gamma^\Delta$ at distance $2$ (in the line graph), and to $U$ corresponds a line $U^*$ which is contained in $\{ V, W \}^{\perp}$. Suppose that $w$ is mapped to a line $\gamma(w) = W'$ in $\mB'$. First assume that $W'$ is not incident with $v'$. Then in $\Gamma'$, $\gamma(V)$ and $\gamma(W)$ are at distance $2$ in the line graph, which means that $W'$ is a point in $\mP'$, contradiction. Next suppose that $W'$ is incident with $v'$. Then one immediately obtains a contradiction, as $(v,U)$ and $(w,U)$ are flags in $\Gamma$, and $\gamma$ preserves incidence (as flags are sent to flags). It follows that $\gamma$ sends collinear points in $\Gamma$ to collinear points in ${\Gamma'}$, and by connectedness, $\gamma$ sends all points in $\mP$ to points in $\mP'$. Since incidence is preserved, all lines go to lines.  So $\gamma$ induces an epimorphism 
\[ \underline{\gamma}: \Gamma \mapsto {\Gamma'}. \]

If on the other hand, $\gamma$ sends some point in $\mP$ to a line in $\mL'$, then connectedness similarly leads to the fact that $\gamma$ sends all points of $\mP$ to lines of $\mL'$ and all lines in $\mL$ to points in $\mP'$. So $\gamma$ induces an epimorphism 
\[ \underline{\gamma}: \Gamma \mapsto {\Gamma'}^D. \] 
(or in other words, $\underline{\gamma}: \Gamma^D \mapsto {\Gamma'}$).  

Conversely, for $n \geq 2$ each epimorphism $\underline{\gamma}$ between generalized $n$-gons defines an epimorphism $\gamma$ between the corresponding doubles. \\

Finally, we have the following result. 

\begin{theorem}
 Let $\gamma: \mS \mapsto \mS'$ be an epimorphism between finite thin generalized $m$-gons of order $(s,1)$ and $(s',1)$, with $m \in \{ 4, 6, 8, 12 \}$.
Then we have the following cases: 
\begin{itemize}
\item[{\rm (a)}] 
if $m = 4$, then each epimorphism $\gamma$ is obtained by doubling any epimorphism between digons $\mD$ and $\mD'$ of order $s$ and $s'$;
\item[{\rm (b)}] 
if $m \in \{ 6, 8, 12 \}$, then either $\gamma$ is an isomorphism (and so $s = s'$), or
$s' = 1$ and all possible epimorphisms $\gamma$ arise from doubling the epimorphisms of the corresponding $3$-gons, $4$-gons, and $6$-gons described in Theorem \ref{GT}, Theorem \ref{JATGQ} and Theorem \ref{JATGH}.
\end{itemize}
\end{theorem}

{\em Proof.}\quad 
 If $\theta$ is a surjection from the point set of a generalized digon $\mD$ onto the point set of a generalized digon $\mD'$ and $\theta'$ is a surjection from the line set of $\mD$ onto the line set of $\mD'$, then these surjections define an epimorphism $\gamma$ from $\mS$ onto $\mS'$
(similarly if we interchange the role of points and lines). The converse is also easy. Part (b) of the theorem is just an accumulation of the discussion prior to the theorem. 
\eop \\

\begin{remark}{\rm
In \cite[4.2.4]{POL} Van Maldeghem describes with coordinates how to construct an epimorphism (and in his specific case, a retraction) from any 
(thick) generalized $n$-gon to some ordinary sub $n$-gon (which has order $(1,1)$).}\\
\end{remark}

\section{Application}
\label{app}

A {\em geometrical hyperplane} in a thick finite generalized quadrangle $\mS$ of order $(s,t)$ is a subgeometry $\mS'$ with the property that each line $U$ of $\mS$ either contains 
one point of $\mS'$, and then $U$ is not a line of $\mS'$, or all its points are points of $\mS'$, and then $U$ is a line of $\mS'$. It is easy to show that there only are 
three types of geometrical hyperplanes in this setting:
\begin{itemize}
\item[\framebox{A}]
$\mS'$ is an ovoid of $\mS$ (and so $\mS'$ does not have lines);
\item[\framebox{B}]
$\mS'$ is the point-line geometry of a set $x^{\perp}$, with $x$ a point of $\mS$ (so the only lines of $\mS'$ are those lines of $\mS$ incident with $x$);
\item[\framebox{C}]
$\mS'$ is a subquadrangle of order $(s,t/s)$. 
\end{itemize}

In \cite{part2}, the following theorem is obtained. 

\begin{theorem}
\label{hyps}
Let $\Delta$ and $\Delta'$ be finite thick generalized quadrangles, and let $\phi: \Delta \mapsto \Delta'$ be a GQ morphism which surjectively maps a geometrical hyperplane $\mG$ of $\Delta$ to 
a geometrical hyperplane $\mG'$ of $\Delta'$:
\begin{equation}
\phi(\mG) = \mG'.
\end{equation}
Then $\phi$ is surjective if $\mG'$ is not of type B, and also if $\mG'$ is not thin if it is of type C. 
\end{theorem}

Now suppose that $\mG'$ is of type C, and thin. Let $\Delta'$ have order $(s',t')$, so that $\mG'$ has order $(s',t'/s')$. Then $s' = t'$, and by Theorem \ref{JATGQ}, it follows that $s' = 1$. This contradicts the assumption.  Hence $\phi$ is always surjective if $\mG'$ is not op type B. \\

\section{What is next?}

The next basic case we wish to handle is that of epimorphisms with source a thick generalized octagon and target a thin generalized octagon. This work will be carried out in a companion paper which is being prepared by the authors. \\



\end{document}